\documentclass[12pt,emtex]{article}

\usepackage{amsfonts}
\usepackage{amssymb}
\usepackage{amscd}
\usepackage{amsmath}

\usepackage[english]{babel}

\newcommand{\BEQ}{\begin{equation}}
\newcommand{\EEQ}{\end{equation}}
\def\bea{\begin{eqnarray}}
\def\eea{\end{eqnarray}}
\def\nn{\nonumber}

\newtheorem{Lem}{Lemma}
\newtheorem{Prop}{Proposition}

\newtheorem{Problem}{Problem}

\newtheorem{Rem}{Remark}

\def\bea{\begin{eqnarray}}
\def\eea{\end{eqnarray}}

\def\bes{\begin{equation*} \begin{split}}
\def\ees{\end{split} \end{equation*}}

\def\C{{\mathcal{ C}}}

\def\A{\mathcal{A}}

\begin{document}
\centerline{\LARGE \bf Deformation quantization } 
\medskip
\centerline{\LARGE \bf of integrable systems}
\vspace{.5cm}
\centerline{G.Sharygin\footnote{gshar@yandex.ru} (MSU\footnote{Lomonosof Moscow State University, GSP-1, Leninskie Gory, Moscow, 119991}, ITEP\footnote{Institute for Theoretical and Experimental Physics,  25, Bolshaya Cheremushkinskaya,
 117218, Moscow} and LDCG\footnote{Laboratory of Discrete and Computational Geometry, Yaroslavl' State University, Sovetskaya st. 14, Yaroslavl', Russia 150000 }) and D.Talalaev\footnote{dtalalaev@yandex.ru} (MSU, ITEP)}
\vskip 1cm
\begin{abstract}
In this paper we address the following question: is it always possible to choose a deformation quantization of a Poisson algebra $\A$ so that certain Poisson-commutative subalgebra $\C$ in it remains commutative? We define a series of cohomological obstructions to this, that take values in the Hochschild cohomology of $\C$ with coefficients in $\A$. In some particular case of the pair $(\A,\C)$ we reduce these classes to the classes of the Poisson relative cohomology of the Hochschild cohomology. We show, that in the case, when the algebra $\C$ is polynomial, these obstructions coincide with the previously known ones, those which were defined by Garay and van Straten (see \cite{Stra}).
\end{abstract}

\section{Introduction}
\subsection{Setting of the problem}
In the theory of integrable systems one starts with a Poisson manifold $M,\,\pi$ (where bivector $\pi$ verifies the equation $[\pi,\,\pi]=0$ for the Schouten-Nijenhuis brackets $[,]$). Given such data, one can introduce the Poisson brackets on the algebra of smooth functions by the rule $\{f,\,g\}=\pi(df,\,dg)$ so that for any Hamiltonian $H\in C^\infty(M)$ we have the dynamics on $M$, given by the formulas $\dot{f}=\{H,\,f\}$ for any $f\in C^\infty(M)$ (in more geometric terms we consider dynamics, determined by the vector field $\pi^{ij}\frac{\partial H}{\partial x^j}$).

In particular, if $\pi$ has maximal rank everywhere, i.e. if the manifold $M$ is symplectic (hence it has even dimension $2n$), then we can use the Liouville theorem: in order to describe the trajectories of a dynamical system, it is enough to find $n$ functionally-independent functions $f_1=H,\,f_2,\dots,\,f_n$, such that $\{f_i,f_j\}=0$. If this is the case, one says, that $f_1,\dots,f_n$ is an integrable system. Generalizing a little, we shall say, that an integrable system on a Poisson manifold $M$ is any algebra $\C\subseteq C^\infty(M)=\A$, such that $\{f,\,g\}=0$ for all $f,\,g\in\C$.

On the other hand, for any Poisson manifold one can define the \textit{deformation quantization} of its algebra of functions (see section \ref{refdef}), a noncommutative algebra $(\A[[\hbar]],\,\ast)$ closely related to the Poisson structure on $C^\infty(M)=\A$. One can say, that an integrable system on $M$ remains integrable after quantization (or determines a \textit{quantum integrable system}), if the subspace $\C[[\hbar]]$ is a commutative subalgebra of $(\A[[\hbar]],\,\ast)$. Observe, that the deformation-quantization approach is not the only one used to define quantum integrable systems. The general theory of quantum integrable systems is a well-developed branch of modern mathematical Physics, we outline its ideas and results in section \ref{refquan}.

Our principal aim in this paper is to find out, if there always exists a quantization of an integrable system in the deformation quantization framework, i.e. the deformation of a Poisson algebra, which preserves its commutative subalgebra, and if it exists to classify all such quantizations. More precisely, let $(\A,\,\{,\})$ be a commutative Poisson algebra, i.e. an algebra with Poisson bracket, verifying the Leibniz rule, and let  $\C$ be its Poisson-commutative subalgebra (i.e. $\C$ is a subalgebra of $\A$ such that the restriction of the Poisson bracket $\{,\}$ on $\C$ vanishes), then we are interested in such a $\ast$-product in $\A[[\hbar]]$  that $\C$ remains commutative subalgebra of $\A[[\hbar]]$ with respect to a $\ast$-product. 

In particular, we can (and shall) assume that the product in $\C$ and in its image in $\A[[\hbar]]$ coincide: although in general, one can ask about the commutativity of $\C[[\hbar]]$ inside $\A[[\hbar]]$, under very mild assumptions (for instance, when the Hochschild cohomology of $\C$ verify the Hochschild-Kostant-Rosenberg theorem), these two conditions are equivalent, up to an isomorphism (see section \ref{firstcalcs}). 

Algebraically this can be written as the following three conditions on an element of the Hochschild complex (see section \ref{hochcoh}) $\Pi\in CH^2_\hbar(\A[[\hbar]])$,
\begin{enumerate}
\item $MC(\Pi)=0$
\item $i(\Pi)=0$
\item $\Pi=\pi +$ higher terms
\end{enumerate}
Here $MC$ denotes the Maurer-Cartan equation
$$
d\Pi+[\Pi,\,\Pi]=0,
$$
where $[,]$ is the Gerstenhaber bracket (see section \ref{hochcoh}). Here we denoted by $i$ the inclusion map $\C\rightarrow \A$ and $i^*$ is the natural extension of $i$ 
$$
i^*:CH^*(\A)\longrightarrow CH^*(\C,\A)
$$
i.e. $i^*$  is the map restricting the polylinear maps from the Hochschild complex of $\A$ to $\C$.

In order to answer the question, if such deformation quantizations exist, we consider the corresponding relative Hochschild complex and define obstructions to such a quantization. We first phrase our results in terms of certain conditions on some cohomology classes, and later rephrase them in terms of the elements of Poisson cohomology on the space of Hochschild cohomology.  We also compare our results with the analogous classes, defined in \cite{Stra}, which turn out to be equal to ours in a simpler situation of a symplectic space with canonical Darboux coordinates. Still another approach to a similar question is contained in a recent paper \cite{Freg} devoted to a wide class of deformation problems of pairs.

The rest of the paper is organized as follows: in sections \ref{refdef} and \ref{refquan} we recall the history, simple facts and notions of deformation quantization and of the theory of quantum integrable systems. In section \ref{hochcoh} we recall the definition of Hochschild cohomology and calculate this cohomology in a particular case. In sections \ref{loter} and \ref{hiter} we define the obstructions in terms of Hochschild cohomology. Finally, in the section \ref{refstra} we describe the relation of these conditions to the results of Garay and van Straten, in particular, we reformulate our results in terms of the Hochschild-Poisson cohomology.
\\

\medskip
{\bf Acknowledgments:} We are grateful to professor Dito for his remarks concerning the draft version of the paper. The work of G.Sh. was supported by the Federal Agency for Science and Innovations of Russian Federation under contract 14.740.11.0081 and by the Russian government project 11.G34.31.0053. The work of D. T. was partially supported by the RFBR grant 11-01-00197-a, by the RFBR-CNRS grant 11-01-93105-a, by the grant of the Dynasty foundation, by the RFBR-OFI-I2 13-01-12401 and by the grant for support of the scientific school 4995-2012.1.

\subsection{Remarks on deformation quantization}
\label{refdef}
The idea of deformation quantization can be traced back to the works of the founders of quantum mechanics: one may argue, that the notion of semi-classical limit of quantum systems describes the latter as a (proto-)deformation of the classical case. Another way to derive the deformation quantization is to consider the Hermann Weyl's quantization formula:
$$
u\mapsto O_u=\int_{\mathbb R^{2l}}\Phi(u)(\xi,\,\eta)\exp(i(\xi^j \partial_j+\eta_k q_k)/\hbar)d^l\xi d^l\eta,
$$
which expresses the operator on $L^2(\mathbb R^l)$, associated to $u\in C^\infty(\mathbb R^{2l})$ in terms of an integral, where $\Phi(-)$ is the inverse Fourier transform, $\partial_j=i\hbar\frac{\partial}{\partial x^j}$ and $q^k$ is the multiplication by $x^k$. The function $u$ here can be interpreted as the symbol of the differential operator $O_u$. The opposite question, how to find an expression of the classical function-symbol of an operator, lead Jose Moyal in 1949 to his famous formula, which expresses the symbol of the product of two operators in terms of the symbols of the factors, which is now called \textit{the Moyal star-product} (at least, so this formula is credited nowadays, although at that time there definitely were other people working on the same subject):
$$
u\ast w= fg+\frac{i\hbar}{2}\pi^{ij}\partial_if\partial_jg+\frac{(i\hbar)^2}{8}\pi^{ij}\pi^{kl}\partial_k\partial_if\partial_l\partial_jg+\dots=m\circ\exp{(\frac{i\hbar}{2}\pi)}(f\otimes g),
$$
where $\pi=\pi^{ij}\partial_i\wedge\partial_j$ is a constant Poisson bivector on $\mathbb R^{2l}$ and $m$ denotes the product of the functions. This formula appeared in 1940-ies and it took some time before it attracted attention of mathematicians. 

The other source of the deformation quantization ideas is its name-sake: deformation theory of complex varieties and its algebraic version, developed in 1960-ies. This theory describes possible ways to pass "continuously" from one algebra, group or some other mathematical object to another. In the framework of this approach various mathematical tools were developed, such as Hochschild homology and cohomology, Gerstenhaber brackets, etc. However, for about two decades this theory has not been applied to the quantum mechanics.

It was probably in the works by M.~Flato and coauthors in mid 1970-ies, see \cite{BFFLS1} and \cite{BFFLS2}, where the physical and mathematical approaches were first synthesized and the following question, which is now generally referred to as \textit{the deformation quantization problem} was formulated:
\begin{Problem}
Let $(M,\,\pi)$ be a Poisson manifold ($\pi$ is Poisson bivector). Find a way to deform formally the product in $C^\infty(M)$, i.e. introduce a new associative product in the space of formal power series $C^\infty(M)[[\hbar]]$, such that it coincides with the original one up to the $\hbar$-terms and the commutator of any two functions $f,g\in C^\infty(M)$ with respect to this product is equal to their Poisson bracket up to $\hbar^2$. Classify such products for a given Poisson structure.
\end{Problem}
One readily sees, that Moyal product gives an example of such noncommutative multiplication on $\mathbb R^{2l}$.

This question alongside with the closely related quantum groups theory (in which one is to find a deformation of the group structure) has been extensively studied in 1980-ies and 1990-ies. Many approaches to it has been developed by various mathematicians: all the machinery of the Hochschild homology, homological algebra, category theory, microlocal analysis and ideas from many other fields were applied. The notable results of this investigations are the Drinfeld's constructions in quantum groups \cite{Dr1}, De Wilde and Lecomte quantization, Fedosov's deformation quantizations of the symplectic manifolds, and  the Kontsevich's quantization theorem \cite{Kont} (one should pay attention both to the original Kontsevich's proof, which amounts to a direct computation by a given euristic formula and to the Tamarkin's proof, based on a general operadic approach).

\subsection{Remarks on quantum integrable models}
\label{refquan}
The theory of quantum integrable models counts numerous examples originated in mathematical physics, namely in spin chains, in condensed matter models, in statistical mechanics  such as Heisenberg magnet, Gaudin system, quantum nonlinear Schr\"odinger equation and many others. More general concept of quantum integrability concerns a pair $\C\subset \A$ of associative algebras where $\C$ is 
 commutative and maximal in an appropriate sense. Due to the algebraic definition there is a deep and fruitful relation of this domain with the representation theory and algebra in general: algebra provides examples of quantum integrable models, and vice versa the methods of quantum integrable models give results in representation theory. 

The main method in this domain is the quantum inverse scattering method (QISM) established in the 70-ies  years of the 20th century by the school of L. D. Faddeev \cite{QISM}. QISM is deeply related with the theory of quantum groups introduced by Drinfeld  \cite{Dr1}. The letter presents a deformation of a classical group in the category of Hopf algebras (there are some generalizations: quasi-Hopf algebras, bialgebras etc.) which is "perpendicular" to the deformation problem of the present paper. Briefly speaking, QISM allows one to construct an integrable system starting with a solution of some structural equation like Yang-Baxter equation related with the corresponding quantum group. 

There is an alternative approach to quantum integrable models which is efficient for a class of models of the Gaudin type. This is a quantum spectral curve method \cite{T} whose principal idea is to consider the quantum integrable model as a deformation of a classical one preserving some additional structures, such as the spectral curve and separated variables. This approach has important advantages against QISM in the solution aspect. The main problem of this work is in a sense analogous: we explore the formal deformations of a classical integrable model up to an equivalence. Besides the models of physical interest we do not discuss the representation of the underlying algebra. Such a difference provides the important distinction in physical properties and will be the subject of future refinement of our approach.

\section{Hochschild complex and deformations}
\subsection{Definitions}
\label{hochcoh}
We have to recall principal facts about the Hochschild complex. Let $\A$ be an associative algebra over a characteristic zero field $\Bbbk$; consider the Hochschild cohomology complex $CH^i(\A)=Hom_\Bbbk(\A^{\otimes i},\A)$ (if $\A$ is a (sub)algebra of the algebra of smooth functions on a manifold, as it is below, we will usually restrict the notion of linear maps to the "local" ones, see remarks in the end of this section). This complex has
\begin{itemize} 
\item a differential $d:CH^i(\A)\rightarrow CH^{i+1}(\A)$
defined on $\varphi:\A^{\otimes i}\rightarrow \A$ as follows
\bea
d\varphi(f_1,\ldots,f_{i+1})&=& f_1 \varphi (f_2, \ldots , f_{i+1})+\sum_{j=1}^{i}(-1)^{j} \varphi(f_1,\ldots,f_j f_{j+1},\ldots, f_{i+1})\nn\\
&+&(-1)^{i+1} \varphi(f_1, . . . , f_{i})f_{i+1}.\nn
\eea
\item the cup-product $\cup:CH^i(\A)\otimes CH^j(\A)\rightarrow CH^{i+j}(\A)$ defined by the formula:
\bea
(\varphi\cup \psi)(f_1,\ldots,f_{i+j})=(-1)^{ij}\varphi(f_1,\ldots,f_i) \psi(f_{i+1},\ldots,f_{i+j})\nn
\eea
\item the Gerstenhaber bracket $[~,~]:CH^i(\A)\otimes CH^j(\A)\rightarrow CH^{i+j-1}(\A)$
\bea
[\varphi,\psi]=\varphi \circ  \psi - (-1)^{(i-1)(j-1)}\psi \circ  \varphi \nn
\eea
where
\bea
(\varphi \circ \psi)(f_1,\ldots,f_{i+j-1})=\sum_{l=1}^{i-1}(-1)^{l(j-1)}\varphi(f_1,\ldots,f_l,\psi(f_{l+1},\ldots,f_{l+j}),\ldots,f_{i+j-1})\nn
\eea
\end{itemize}
The bracket with the differential make the Hochschild complex a differential graded Lie algebra, while the differential and the cup product together define the noncommutative differential graded algebra. Moreover being restricted to cohomology the cup-product and the bracket provide a structure of Gerstenhaber algebra on $HH^*(\A)$. Another important fact about the Hochschild cohomology is that when $\A$ is an algebra of smooth functions on a manifold, its Hochschild cohomology as a Gerstenhaber algebra can be described in classical terms: according to the well-known Hochschild-Kostant-Rosenberg theorem (see \cite{Loday} for example) it is equal to the algebra of polyvector fields on the manifold with the bracket given by the Schouten-Nijenhuis bracket.

Hochschild complex plays a prominent role in the deformation problem: one can regard the deformed multiplication in $\A[[\hbar]]$ as a formal series
$$
a*b=ab+\hbar B_1(a,\,b)+\hbar B_2(a,\,b)+\dots.
$$
Then the associativity condition for $\ast$ can be expressed as the Maurer-Cartan equation on the element $\Pi=\hbar B_1+\hbar B_2+\dots$ in the $\hbar$\/-linear Hochschild complex of $\A[[\hbar]]$, i.e. $MC(\Pi)=0$ (see the introduction).

In what follows we shall assume, that $\A$ is an algebra of functions on a Poisson manifold $M$, $\C$ its Poisson-commutative subalgebra (i.e. $\{f,\,g\}=0$ for all $f,\,g\in\C$). For instance, we can take $\C=\rho^*(C^\infty(X))$ for a map $\rho:M\to X$, intertwining the given Poisson structure on $M$ with the trivial structure on $X$. Or else $\C$ can be the algebra of integrals of an integrable system (in particular, this is especially interesting if $M$ is a symplectic manifold). 

Throughout the paper \textit{we shall consider "local" (with respect to $M$) Hochschild complex}, i.e. the complex consisting of such cochains $\varphi:\A^{\otimes n}\to \A$, that $\varphi(f_1,\dots,f_n)(x)=\varphi(g_1,\dots,g_n)(x)$ if there exists an open neighborhood $U$ of $x$, in which $f_i\equiv g_i,\,i=1,\dots,n$. One can show, that (on smooth manifolds) this is equivalent to the condition, that all the cochains in $CH^*(\A)$ are given by the polydifferential operatos on $M$. In particular, even when we speak about cochains on $\C$, we assume that they are local on $M$; in the terms of differential operators, one can formulate this as follows: every polydifferential operator on the sub algebra $\C$ with values in $\A$ extends to a polydifferential operator on whole $\A$ (i.e. we need a bit more than a linear extension of operators). This condition is evidently fulfilled, when $\C$ is given by the inverse image of a projection, $\C=\rho^*(C^\infty(X))$.

\subsection{A variant of the HKR theorem}
We use the assumptions and notation from previous section. Consider the following exact sequence of Hochschild complexes:
\begin{equation}
\label{eqexseq}
0\rightarrow IQ^*(\A,\,\C)\rightarrow CH^*(\A)\stackrel{i}{\rightarrow} CH^*(\C,\A)\rightarrow 0.
\end{equation}
Here $IQ^*(\A,\,\C)$ denotes the kernel of the natural restriction map $i$; it can be described as the set of all cochains $\varphi\in CH^*(\A)$ that vanish if all the arguments are from $\C$. The exactness of this sequence on the right is evident on the level of the conventional Hochschild complex, i.e. when the elements are just polylinear maps, on the level of polydifferential operators it was conjectured in the previous section. We are going to describe the corresponding cohomological long exact sequence in the case, when $\A=C^\infty(M)$ and $\C=C^\infty(X)$ for a smooth submersion 
$$
\rho:M\rightarrow X.
$$
To this end consider the exact sequence of vector bundles
$$
0\to T^{vert}_\rho M\to TM\to T^{hor}_\rho M\to 0,
$$
where $T^{vert}_\rho M$ is the kernel of the differential of $\rho$ and $T^{hor}_\rho M=TX/T^{vert}_\rho M$, which we can identify with the pullback $\rho^*TX$.
\begin{Prop}
\label{prop1}
Cohomology of the complexes that appear in \eqref{eqexseq} are given by the formulas
$$
HH^*(\A)\cong\wedge^*TM,\ H^*(\C,\A)\cong\wedge^*T^{hor}_\rho M,\ H^*(IQ^*(\A,\,\C))\cong\langle T^{vert}_\rho M\rangle,
$$
where $\wedge^*TM$ (resp. $\wedge^*T^{hor}_\rho M$) denotes the algebra of polyvector fields on $M$ (resp. the algebra of "horizontal" polyvector fields on $M$, which can be regarded as the pullbacks of polyvector fields on $X$), and $\langle T^{vert}_\rho M\rangle$ denotes the kernel $\ker \wedge(d\rho)$, the ideal in $\wedge^*TM$, generated by $T^{vert}_\rho M$. The long exact sequence of cohomology, associated with \eqref{eqexseq} splits into short exact sequences of the form
$$
0\to \langle T^{vert}_\rho M\rangle^k\to \wedge^kTM\to \wedge^kT^{hor}_\rho M\to 0. 
$$
\end{Prop}
\noindent\textit{Proof}
The isomorphism $HH^*(\A)=\wedge^*TM$ is the conclusion of the Hochschild-Kostant-Rosenberg theorem. Further, since our complexes are local in $M$, we can restrict the exact sequence to any open neighborhood in $M$ and use partition of unity to restore the general result from the local ones (\`a la Mayer-Vietoris sequence, see, for instance \cite{Bry89}). Thus we can assume that $M=X\times F$ for a fibre $F$. Then the equality  $H^*(\C,\A)\cong\wedge^*T^{hor}_\rho M$ becomes quite evident: the Hochschild-Kostant-Rosenberg map
$$
\begin{aligned}
\chi_{HKR}&:\wedge^*TM\to CH^*(\A),\\ 
\chi&(X_1\wedge\dots\wedge X_n)(f_0,\dots,f_n)=\frac{1}{n!}\sum_{\sigma\in S_n}(-1)^\sigma f_0X_1(f_{\sigma(1)})\dots X_n(f_{\sigma(n)}),
\end{aligned}
$$
induces a map $\chi_{HKR}':\wedge^*T^{hor}_\rho M\to CH^*(\C,\,\A)$, which commutes with the differential and is clearly an isomorphism. The rest follows from the exactness of the long sequence in Hochschild cohomology.
\hfill$\square$

Here is a couple of important observations, that follow from this proposition:
\begin{enumerate}
\item All the maps 
$$
HH^k(\A)\to H^k(\C,\A)
$$
are epimorphic;
\item The proposition stays true, if instead of the submersion $\rho$ we have only an integrable distribution $\omega$, so that $T^{vert}_\omega M$ consists of vectors in $\omega$ and $T^{hor}_\omega M=TM/T^{vert}_\omega M$ and take $\C$ to be the algebra of functions on $M$, eliminated by the vertical vector fields.
\end{enumerate}

\section{Obstructions and calculations}
\subsection{Deformation problem}
\label{firstcalcs}
Let $\ast$ be a deformed product (for example a product, given by Kontsevich's theorem) on a Poisson algebra $\A$, in other words it is a deformation of a Poisson algebra $\A$. From the point of view of the Hochschild complex, this is an element $\Pi$ in $CH^*(\A[[\hbar]])$, that verifies the Maurer-Cartan equation and (since we need to keep trace of the Poisson structure) begins with the Poisson bracket. Also recall, that two deformations $\ast_1$ and $\ast_2$ are called \textit{equivalent}, if there exists a formal power series of differential operators $D=id+\hbar D_1+\hbar^2D_2+\dots$, such that  
$$
D(a\ast_2 b)=D(a)\ast_1 D(b).
$$
Due to Kontsevich's formality theorem \cite{Kont} $\ast$\/-products of this sort always exist and their equivalence classes coincide with the equivalence classes of formal Poisson structures under similar formal isomorphisms. Recall, that a formal Poisson structure is the formal power series of bivectors: 
$$
\Pi=\pi_0+\hbar\pi_1+\hbar^2\pi_2+\dots,
$$
that verify the usual properties of the Poisson bivectors in the space of formal series.

So we ask, if a given class of $\ast$\/-products contains a product, trivial on $\C$, in other words, given a $\ast$\/-product, we are going to look for an ($\hbar$\/-linear) automorphism $D$ of the space $\A[[\hbar]]$ such that 
\bea
D^{-1}(D(a)*D(b))=a b\ \forall a,b\in\C\nn
\eea
or
\bea
\label{mor}
D(a)*D(b)=D(a b)\qquad \forall a,b\in\C,
\eea
and do this term by term in expansion on $\hbar.$

Let us introduce some notations for the automorphism and the deformation series:
\bea
D(a)&=&a+\hbar D_1(a)+\hbar^2 D_2(a)+\ldots \nn\\
a\ast b&=&a b+\hbar B_1(a,b)+\hbar^2 B_2(a,b)+\ldots \nn
\eea
Observe, that $B_1(a,b)$ is not necessarily equal to the Poisson bracket $\{,\}$. However, we can always find a formal diffeomorphism $D$, which will intertwine any given $\ast$\/-product with a $\ast$\/-product in which $B_1(a,b)=\{a,b\}$, so we shall fix this equality from now on. In effect, in a similar way one can show, that every commutative $\ast$\/-product on $\C[[\hbar]]$ can be replaced by the trivial one $a\ast b=a b$ in $\C$.

\subsection{$\hbar^2$-term}
\label{loter}
Expanding both sides of  (\ref{mor}) and collecting terms at $\hbar$ and $\hbar^2$ one obtains
\bea
\hbar:&& aD_1(b)+D_1(a)b=D_1(ab)\nn\\
\hbar^2:&& B_2(a,b)+B_1(a,D_1(b))+B_1(D_1(a),b)+D_1(a)D_1(b)+aD_2(b)+D_2(a)b=D_2(ab)\nn
\eea
The first equality means that $D_1$ is a derivation on a subalgebra $\C$ with values in $\A$. Denoting by $d$ the Hochschild differential we reduce the second equality to the following one:
\bea
\label{prep}
B_2(a,b)=d D_2(a,b)-D_1\cup D_1 (a,b)-[D_1,B_1](a,b)
\eea
here the bracket means the Gerstenhaber bracket on the Hochschild complex and $\cup$ is the cup-product. Let us remember that this equality should fulfill for all $a,\,b\in\C$. Let us also emphasize that the second term may be nontrivial despite  the fact that $B_1|_\C=0$.

Let us also recall some consequences from the associativity of the $\ast$-product for first terms of the deformation series:
\bea
dB_1(a,b,c)&=&0;\nn\\
dB_2(a,b,c)&=&B_1(B_1(a,b),c)-B_1(a,B_1(b,c)).\nn
\eea
The first one means that the element $B_1$ is closed in Hochschild complex, in particular this equation holds, if one chooses the Poisson bivector $\pi$ for the $B_1$, (as we have done before). Moreover $B_1(a,b)=\{a,b\}=0$ for all $a,\,b\in\C$ hence $dB_2=0$ when restricting to $\C$. Hence \textit{the second term in the deformation series for the $\ast$-product, when restricted to $\C$, defines an element $i^*(B_2)$ in the cohomology space $H^2(\C,\,\A)$}, which we shall denote by $B_2$, when it does not cause ambiguity.

{\Lem \label{lem1} The necessary condition for the automorphism $D$ with property \ref{mor} to exist
is that $B_2$ lies in the sum of images of $\cup:H^1(\C,\,\A)\otimes H^1(\C,\,\A)\to H^2(\C,\,\A)$ and the Gerstenhaber action $[\pi,\,\cdot]:H^1(\C,\,\A)\to H^2(\C,\,\A)$.}

Below we shall say more on the definition of the map $[\pi,\,\cdot]$, see section \ref{refstra}. Now one can regard this term as the restriction of the Gerstenhaber bracket.

Evidently this condition is not sufficient for the full problem but this shows that taking the two-term automorphism 
$$
D=a+\hbar D_1+\hbar^2 D_2,
$$ 
whose coefficients  $D_1,\,D_2$ satisfy \eqref{prep} one can find a deformation of  the $\ast$-product such that $B_2(a,b)=0$ for all $a,b\in\C$. 
{\Rem $ $\\

\vspace{-.5cm}
\begin{enumerate}
\item In the conditions of the Proposition \ref{prop1} the square of any element $[D]\in H^1(\C,\A)$ is equal to $0,$ hence the necessary condition in lemma \ref{lem1} is that $i(B_2)$ is $[\pi,\,\cdot]$-exact in $H^2(\C,\,\A)$.
\item The condition of the lemma \ref{lem1} is in effect necessary and sufficient for the elimination of the first two terms of the $\ast$\/-product on $\C$, unlike the conditions that we shall give in the following section, which are merely sufficient.
\end{enumerate}
}

\subsection{Higher terms}
\label{hiter}
Suppose, the obstruction we found in the previous section is trivial, in particular, we can find an $\ast$\/-product so that
$$
B_2(a,b)=0,\ \forall a,b\in\C.
$$
Let us consider an automorphism $D$ with $D_1=0$ and $D_2,\,D_3,\dots$ such that \eqref{mor} fulfills up to $\hbar^3.$ Then the associativity condition implies
$$
dB_3(a,b,c)= [B_2,B_1](a,b,c)
$$
which is $0$ if $a,\,b,\,c\in\C.$ That is $B_3$ is closed in $CH^2(\C,\,\A)$.

Substituting $D$ and collecting terms with different powers of $\hbar$ we obtain
$$
\begin{aligned}
dD_2(a,b)&=0,\\
B_3(a,b)&=d D_3(a,b)-[D_2,B_1](a,b),
\end{aligned}
$$
for all $a,b\in\C$. Hence there exists a deformation of $\ast$\/-product by a formal diffeomorphism, such that \eqref{prep} fulfills up to $\hbar^4$, if the $H^2(\C,\A)$-cohomology class of $B_3$ lies in the image of the Gerstenhaber bracket $[\pi\, ,\cdot\,]:H^1(\C,\,\A)\to H^2(\C,\,\A)$.

This observation can be generalized.
{\Lem
Let $\ast_n$ be a $\ast$-product equivalent to $\ast$ such that
$$
B_1(a,b)=B_2(a,b)=\dots=B_n(a,b)=0,\ \forall a,b\in\C.
$$
then, 
\begin{itemize}
\item the term $B_{n+1}:\C\otimes\C\to\A$ is a Hochschild cocycle;
\item the $\ast$-product $\ast_n$ can be deformed to $\ast_{n+1}$ if the class of  $B_{n+1}$ in $H^2(\C,\,\A)$ lies in the image of $[\pi\, ,\cdot\,]:H^1(\C,\,\A)\to H^2(\C,\,\A)$.
\end{itemize}
}
{\Rem
We should once again warn the reader, that the conditions of the present section are sufficient, but not in general necessary. In fact, in order to get necessary and sufficient conditions, we should have dropped the assumption that the formal diffeomorphism begins with zeros. This would have caused additional terms in the calculations. For instance, if we do so for the term $B_3$, we would obtain the following set of relations:
\begin{align*}
dD_1&=0,\\
dD_2&=D_1\cup D_1+ [\pi,\,D_1],\\
\intertext{and the following relation for $B_3$ (we use the notation $B_1$ instead of $\pi$ for the sake of uniformity)}
B_3&=[D_1,B_2]+[D_2,B_1]+D_1\cup D_2+D_2\cup D_1+dD_3.
\end{align*}
Similar formulas can be written for arbitrary level $n$. As one can see, the necessary condition would involve the existence of a very special formal diffeomorphism, which would kill the given class in $H^2(\C,\,\A)$. It is intriguing to find a good geometric interpretation of the corresponding relations, but for a time being we shall restrict our attention to the sufficient conditions we have obtained. 
}

\subsection{Poisson cohomology}
\label{refstra}
In previous sections we have defined a series of cohomology classes $B_2,\,B_3,\dots\in H^2(\C,\,\A)$, which belong to the image of the Gerstenhaber bracket $[\pi\, ,\cdot\,]$, if the integrable system can be quantized. Now we are going to describe the same condition in a bit more intrinsic way.

Recall the definition of Poisson cohomology: let $(M,\,\pi)$ be a Poisson manifold, $\wedge^*TM$ will denote the space of polyvector fields on $M$. One can define the differential $d_\pi$ of degree $+1$ on this space by equation
$$
d_\pi:\wedge^kTM\to\wedge^{k+1}TM,\ d_\pi(T)=[\pi,\,P],
$$
where $P$ is a polyvector field and $[~,~]$ is the Schouten-Nijenhuis bracket. The equality $d_\pi^2=0$ follows from the Jacobi identity which is equivalent to the equation $[\pi,\,\pi]=0$, which determines the Poisson bivector. The cohomology of this complex is called the Poisson cohomology of $(M,\,\pi)$. It is closely related to the Poisson homology of Brylinski, e.g. see \cite{Cic}.

Let now $V\subseteq TM$ be a distribution of subspaces in $TM$ such that for any vector field $X\in V$ we have $[\pi,\,X]\in [V]$ where $[V]$ is defined as the kernel of the exact sequence:
\bea
0\rightarrow [V]\rightarrow \wedge^* TM\rightarrow \wedge^* \left(TM/V\right)\rightarrow 0.\nn
\eea  
Then we can define a version of the Poisson differential both on the space $[V]$ and on the space of sections of the quotient-bundle $H=TM/V$. Indeed, for any vector field $Y\in H=TM/V$, we choose a representative $\tilde Y\in TM$ and define $d_\pi Y=[\pi,\,\tilde Y](\mathrm{mod}\,V)$. From assumptions we mae it follows at once that the result doesn't depend on the choice of representative. Since the Schouten bracket on higher dimensional polyvector fields is defined with the help of the Leibniz rule, we obtain the differential $d^H_\pi$ on $\wedge^*H$. Similarly, the same formula defines a differential on the kernel of the projection $\wedge^* TM\to \wedge^* H$, which we shall also denote by $d^V_\pi$. More generally, if $V$ is an integrable distribution, then we can define an action of the Schouten-Poisson algebra of polyvector fields with values in $V$ on the space $\wedge^*H$: the same consideration shows, that the usual formulas give a well-defined result. 

We shall call the cohomology of $(\wedge^*H,\,d_\pi^H)$ the \textit{relative Poisson cohomology} of $M$ modulo $V$. In particular, in the case we considered in the section \ref{hochcoh}, we showed that the Hochschild cohomology of the pair $\A,\,\C$ is equal to $\wedge^*T^{hor}_\rho M$, where $T^{hor}_\rho M$ is the quotient-bundle of $TM$ modulo a distribution, induced by a projection (or, more generally, modulo any integrable distribution, see remark following the proof of the proposition \ref{prop1}). In this case we shall denote the differentials $d^V_\pi$ and $d^H_\pi$ by $d^{vert}_\pi$ and $d^{hor}_\pi$ respectively. Recall now that the image of the Gerstenhaber bracket under the identification of the Hochschild-Kostant-Rosenberg theorem is the Schouten bracket of polyvector fields.

The following proposition is in certain sense an algebraic analogue of the remarks, concerning the differential $d_\pi$:
\begin{Prop}
Let $\A$ be the algebra of smooth functions on a manifold and $\C$ its subalgebra, defined as in the conditions of proposition \ref{prop1}. Then the Gerstenhaber bracket in $CH^*(\A)$ can be restricted to the subcomplex $IQ^*(\A,\,\C)$ and for any $\varphi\in IQ^p(\A,\,\C)$ and $\xi\in CH^q(\C,\,\A)$ the formula
$$
[\varphi,\,\xi](a_1,\dots,a_{p+q-1})=\sum_{i=1}^p(-1)^{iq+1}\varphi(a_1,\dots,a_{i-1},\xi(a_i,\dots,a_{i+q-1}),a_{i+q}\dots,a_{p+q-1})
$$
determines an action of the Lie algebra $IQ^*(\A,\,\C)$ on $CH^*(\C,\,\A)$. The image of this action on $\wedge^*T^{hor}_\rho M$ is given by the Schouten bracket on polyvector fields.
\end{Prop}
\noindent\textit{Proof}
It is enough to observe, that in the case we consider the second term of the usual Gerstenhaber bracket of $\varphi$ and $\xi$ (or rather the restriction of $\varphi$ to $\C$) should vanish, since $\varphi\in IQ^*(\A,\,\C)$. The rest are the classical results of Gerstenhaber.
\hfill$\square$

In particular, the Gerstenhaber bracket with $\pi\in IQ^*(\A,\,\C)$ in the view of the results of proposition \ref{prop1} induces the differentials $d_\pi,\,d_\pi^V$ and $d_\pi^H$. Now the conclusions of our previous sections can be reformulated as follows:
\begin{Prop}
Consider an integrable system $(\A,\,\C,\,\{,\})$, where $\A$ and $\C$ are as in the conditions of proposition \ref{prop1}. Then the obstruction classes $B_n\in H^2(\C,\,\A)=\wedge^2T^{hor}_\rho M$ are closed with respect to $d^{hor}_\pi$ and the deformation of integrable system exists if they are exact. 
\end{Prop}
\noindent\textit{Proof}
The only thing that needs checking is the closedness of $B_n$ for all $n$. But this follows from the associativity equation: a direct computation shows that is $B_1=\dots=B_{n-1}=0$ on $\C$, then
$$
[\{,\},\,B_n]=d(B_{n+1}).
$$
\hfill$\square$

In what follows we shall denote the corresponding classes in Poisson cohomology by $\widetilde{B_n}$.

Let now $\C\cong\mathbb R[x_1,\dots,x_n]$ with generators $x_i$ given by functions $f_i\in C^\infty(M)$. One can use the Koszul resolution of $\C$ to compute the Hochschild cohomology. Recall (\cite{Loday}) that this resolution is given by
$$
K^*(\C)=\oplus_{i=0}^n \C\otimes\wedge^i\mathbb R^n\otimes\C,
$$
with differential given by $d(x\otimes v\otimes y)=xv\otimes y-x\otimes vy$ on $\C\otimes\mathbb R^n\otimes\C$ (where we identify $x_i\in\mathbb R^n$ with $f_i\in\C$) and is extended to whole $K^*(\C)$ by the Leibniz rule. It is straightforward to see now, that in this case
$$
H^*(\C,\,\A)\cong\A\otimes\wedge^*\mathbb R^n.
$$
If $\C\cong\mathbb R[x_1,\dots,x_n]$ is Poisson-commutative subalgebra in $\A=C^\infty(M)$ we can consider map $M\to\mathbb R^n$, given by $x\mapsto(f_1(x),\dots,f_n(x))$. This map is submersion if $f_i$ are functionally independent, so that $\C$ can be regarded as the algebra of functions, eliminated by vertical vector fields of a foliation, verifying the conditions of proposition \ref{prop1}. Thus we can consider the differential $d^{hor}_\pi$. It is straightforward to see, that it is given by the formula
\begin{equation}
\label{eqdfpoi}
d^{hor}_\pi(w\otimes v)=\sum_{i=1}^n\{f_i,\,w\}\otimes e_i\wedge v
\end{equation}
for all $w\in\A$, if $e_i\in\mathbb R^n$ are the corresponding basis elements. In effect, for any element $f\otimes e_i\in\A\otimes\mathbb R^n$ we choose a representative vector field $\tilde e_i$ of $e_i$ on $M$ (we can do it locally, assuming, that the support of $f$ is small enough; it is sufficient, since both formulas-definitions of $d_\pi^{hor}$ are local). Now if $\pi=X\wedge Y$ on the chosen subset, where both $X$ and $Y$ are tangent to the fibers of the foliation, we conclude, that the representative $\tilde e_i$ can commute with $X$ and $Y$, so the formula \eqref{eqdfpoi} holds.

Thus, the complex $(\wedge^*T^{hor}_\rho M,\,d^{hor}_\pi)$ in this case coincides with the complex of Garay and van Straten (see \cite{Stra} and definitions therein). It is now easy to prove the following
\begin{Prop}
The classes $\widetilde{B_n}$ we have defined coincide with the classes $\chi_n$ of Garay and van Straten.
\end{Prop}
\noindent\textit{Proof}
In their paper Garay and van Straten deform the series, corresponding to $f_i$ in $\A[[\hbar]]$ so that $[f^{(n-1)}_i,\,f^{(n-1)}_j]=o(\hbar^{n})$ (the superscript $(n-1)$ denotes the $n-1$\/-st stage of the process). In order to obtain these series in our setting, use the deforming series $D(f_i)$ instead of $f_i^{(n)}$ where $D$ is determined by the $n$\/-th step of the iterative process from sections \ref{loter} and \ref{hiter}, then the commutator relations will follow from the condition on deformed multiplication. The classes of Garay and van Straten were given by
$$
\chi_n=\sum_{i,j}[f^{(n-1)}_i,\,f^{(n-1)}_j]^{(n)}e_i\wedge e_j.
$$
Here $[,]^{(n)}$ denotes the coefficient at $\hbar^n$ in the corresponding formula. Now the $n$\/-th degree in $\hbar$ of the commutators of elements $D(f_i)$ and $D(f_j)$ in $\A[[\hbar]]$ is given by $B_n(f_i,\,f_j)$. It is now enough to recall the formula of the map $\chi_{HKR}'$ from proposition \ref{prop1} to obtain the result.
\hfill$\square$

\section{Conclusion}
In conclusion, we would like to discuss some further questions, concerned with the classification of quantum integrable systems, as well as to point out the direction of our further investigations.

First of all, the already classical results of Kontsevich can be reinterpreted in terms of formality statement: in his paper \cite{Kont} he in effect constructs an $L_\infty$\/-quasi-isomorphism between the differential graded Lie algebra of Hochschild cochains (with respect to the Gerstenhaber bracket) and the algebra of polyvector fields on a manifold (with respect to the Schouten bracket). In our case, we have an exact sequence of Hochschild complexes \eqref{eqexseq}, rather than just one complex and the corresponding exact sequence of the cohomology. Kontsevich's theorem shows, that the complex in the middle is formal. Now the problem we address in this paper can be reformulated as the following question about formality of another complex in the exact sequence: observe, that the Poisson structure $\pi$ as a Hochschild cochain belongs to $IQ^2(\A,\,\C)$; thus the deformation problem we consider will be solved if we can prove that $IQ^*(\A,\,\C)$ is formal. In fact, if this is so, then for any formal Poisson structure $\pi\in H^*(IQ^2(\A,\,\C))$, we shall have a formal solution to the Maurer-Cartan equation $\Pi\in IQ^2(\A,\,\C)$, extending it, just like in the Kontsevich's theorem.

It is not quite clear, if the complex $IQ^*(\A,\,\C)$ is formal or not. In the paper of Garay and van Straten (see \cite{Stra}) it was shown that the introduced obstructions vanish on symplectic manifolds, however it is not at all easy to say, whether the same holds in a generic case. In our attempts to clarify it we calculated few first obstructions in Kontsevich's formula in some particular cases, which all turn out to be trivial. One should observe, that the classes we obtain belong to the cohomology of the right-hand complex in the exact sequence, while the formality problem is concerned with the complex on the left. The reason for this might be in the fact, that the exact sequence \eqref{eqexseq} represents an extension in the category of differential Lie algebras, thus the formality might be closely related to the class of this extension in the derived category. However, so far we are only able to suggest some speculations on this, rather intriguing subject.

Another interesting question is to find an explicit formula for the deformation quantization in this case. We want to note that our efforts to construct an explicit formula for the quantization, by analogy with the Kontsevich quantization formula for the Poisson algebra, have not been successful. These questions form a basis for further research.

\end{document}